\documentclass[11pt,a4paper]{article}
\usepackage[latin1]{inputenc}
\usepackage{amsmath}
\usepackage{amsfonts}
\usepackage{amssymb}
\usepackage{theorem}
\usepackage{mathrsfs}

\newtheorem{dfn}{Definition}[section]
\newtheorem{lem}[dfn]{Lemma}
\newtheorem{theor}[dfn]{Theorem}

\newtheorem{rem}[dfn]{Remark}
\newtheorem{cor}[dfn]{Corollary}
\newtheorem{ex}[dfn]{Example}

\newtheorem{pf}[dfn]{Proof.}
\title{ Advances in Startpoint Theory for quasi-pseudometric spaces. 
}
\author{ Ya\'e Ulrich Gaba, \\ gabayae2@gmail.com, \\ Department of Mathematics and Applied Mathematics,\\
University of Cape Town, Rondebosch 7701, South Africa.}

\begin{document}
\maketitle

\begin{abstract}
 This paper presents some startpoint (endpoint, fixed point) theorems for mutli-valued maps that generalize recent  results proved by Y. U. Gaba \cite{rico, ricoo}.
\end{abstract}

\section{Introduction and prelimaries}

\begin{dfn}
Let $X$ be a non empty set. A function $d:X \times X \to [0,\infty)$ is called a \textbf{quasi-pseudometric} on $X$ if:
\begin{enumerate}
\item[i)] $d(x,x)=0 \quad \forall \ x \in X$, 
\item[ii)] $d(x,z) \leq d(x,y) + d(y,z) \quad \forall\  x,y,z \in X $. 
\end{enumerate}
Moreover, if $d(x,y)=0=d(y,x) \Longrightarrow x=y$, then $d$ is said to be a \textbf{$T_0$-quasi-pseudometric}. The latter condition is referred to as the $T_0$-condition.
\end{dfn}

\begin{rem} \hspace*{0.5cm}  
\begin{itemize}
\item Let $d$ be a quasi-pseudometric on $X$, then the map $d^{-1}$
defined by $d^{-1}(x,y)=d(y,x)$ whenever $x,y \in X$ is also a quasi-pseudometric on $X$, called the \textbf{conjugate} of $d$. In the literature, $d^{-1}$ is also denoted $d^t$ or $\bar{d}$.
\item It is easy to verify that the function $d^s$ defined by $d^s:=d\vee d^{-1}$, i.e. $d^s(x,y)=\max \{d(x,y),d(y,x)\}$ defines a \textbf{metric} on $X$ whenever $d$ is a $T_0$-quasi-pseudometric.
\end{itemize}

\end{rem}

Let $(X,d)$ be a quasi-pseudometric space. 

For $x \in X$ and $\varepsilon > 0$, $$C_{d}(x,\varepsilon)=\lbrace y \in X: d(x,y) < \varepsilon \rbrace $$ denotes the open $\varepsilon$-ball at $x$. The collection of all such balls yields a base for the topology $\tau (d)$ induced by $d$ on $X$. Hence, for any $A \in X$, we shall respectively denote by $ int_{\tau (d)}A$ and $cl_{\tau (d)}A$ the interior and the closure of the set $A$ with respect to the topology $\tau (d)$.

Similarly, for $x \in X$ and $\varepsilon \geq 0$, $$C_{d}(x,\varepsilon)=\lbrace y \in X: d(x,y) \leq \varepsilon \rbrace $$ denotes the closed $\varepsilon$-ball at $x$. 
In the case where $(X,d)$ is a $T_0$ quasi-pseudometric space, we know that $d^s$ defined by $d^s:=d\vee d^{-1}$, i.e. $d^s(x,y)=\max \{d(x,y),d(y,x)\}$ defines a metric on $X$. Hence, we shall say that a subset $E\subset X$ is {\it join-closed} if it is {$\tau(d^s)$-closed}, i.e. closed with respect to the topology generated by $d^s$.
The topology $\tau(d^s)$ finer than the topologies $\tau(d)$ and $\tau(d^{-1})$.

\begin{dfn} Let $(X,d)$ be a quasi-pseudometric space.
The convergence of a sequence $(x_n)$ to $x$ with respect to $\tau(d)$, called \textbf{$d$-convergence} or \textbf{left-convergence} and denoted by $x_n \overset{d}{\longrightarrow} x$,   
is defined in the following way
\begin{equation}
x_n \overset{d}{\longrightarrow} x \Longleftrightarrow d(x,x_n) \longrightarrow 0 .
\end{equation}

Similarly, the convergence of a sequence $(x_n)$ to $x$ with respect to $\tau(d^{-1})$, called \textbf{$d^{-1}$-convergence} or \textbf{right-convergence} and denoted by $x_n \overset{d^{-1}}{\longrightarrow} x$,

is defined in the following way
\begin{equation}
x_n \overset{d^{-1}}{\longrightarrow} x \Longleftrightarrow d(x_n,x) \longrightarrow 0 .
\end{equation}

Finally, in a quasi-pseudometric space $(X,d)$, we shall say that a sequence $(x_n)$ \textbf{$d^s$-converges} to $x$ if it is both left and right convergent to $x$, and we denote it as $x_n \overset{d^{s}}{\longrightarrow} x$ or $x_n \longrightarrow x$ when there is no confusion.
Hence
\[
x_n \overset{d^{s}}{\longrightarrow} x \ \Longleftrightarrow \  x_n \overset{d}{\longrightarrow} x \ \text{ and }\ x_n \overset{d^{-1}}{\longrightarrow} x.
\]

\end{dfn}

\begin{dfn}
A sequence $(x_n)$ in a quasi-pseudometric $(X,d)$ is called
\begin{itemize}
\item[(a)] \textbf{left $d$-Cauchy} with respect to $d$ if for every $\epsilon >0$, there exist $x\in X$ and $n_0 \in \mathbb{N}$ such that 
$$ \forall \  n \geq n_0  \quad d(x,x_n)<\epsilon ;$$

\item[(b)] \textbf{left $K$-Cauchy} with respect to $d$ if for every $\epsilon >0$, there exists $n_0 \in \mathbb{N}$ such that 
$$ \forall \  n,k: n_0\leq k \leq n \quad d(x_k,x_n )< \epsilon  ;$$

\item[(c)] \textbf{ $d^s$-Cauchy} if for every $\epsilon >0$, there exists $n_0 \in \mathbb{N}$ such that 
$$ \forall n, k \geq n_0 \quad d(x_n,x_k ) < \epsilon .$$
\end{itemize}

Dually, we define in the same way, \textbf{right $d$-Cauchy}  and \textbf{right $K$-Cauchy} sequences.
\end{dfn}

\begin{rem}\hspace*{0.5cm} 
\begin{itemize}
\item $d^s$-Cauchy $\Longrightarrow$ left $K$-Cauchy $\Longrightarrow$ left $d$-Cauchy.The same implications hold for the corresponding right notions. None of the above implications is reversible.
\item A sequence is left $K$-Cauchy with respect to $d$ if and only if it is right $K$-Cauchy with respect to $d^{-1}$. 
\item A sequence is $d^s$-Cauchy if and only if it is both left and right $K$-Cauchy.

\end{itemize}
\end{rem}

\begin{dfn}
 A quasi-pseudometric space $(X,d)$ is called

\begin{itemize}
\item \textbf{left-$K$-complete} provided that any left $K$-Cauchy sequence is $d$-convergent,
\item {\bf left Smyth sequentially complete} if any left $K$-Cauchy sequence is $d^s$-convergent.
\end{itemize} 
  \end{dfn}
The dual of notions of \textbf{right-completeness} are easily derived from the above. 

\begin{dfn}
 A $T_0$-quasi-pseudometric space $(X,d)$ is called \textbf{bicomplete} provided that the metric $d^s$ on $X$ is complete.
\end{dfn}

\vspace*{0.5cm}
As usual, a subset $A$ of a quasi-pseudometric space $(X,d)$ will be called {\it bounded} provided that there exists a positive real constant $M$ such that $d(x,y)<M$ whenever $x,y\in A$. Note that 
a subset $A$ of $(X,d)$ is bounded if and only if there are $x\in X$ and $r,s\geq 0$ such that $A\subset C_{d}(x,r)\cap C_{d^{-1}}(x,s)$.

\vspace*{0.5cm}

We set $\mathscr{P}_0(X):=2^X \setminus \{ \emptyset\}$ where $2^X$ denotes the power set of $X$. For $x\in X$ and $A,B \in \mathscr{P}_0(X)$, we set:

$$ d(x,A)= \inf\{ d(x,a),a\in A\}  \text{  \ \
 and  \ \  }  d(A,x)= \inf\{ d(a,x),a\in A\},$$
and define $H(A,B)$ by 

$$H(A,B)= \max \left\lbrace \underset{a\in A}{\sup}\ d(a,B), \underset{b\in B}{\sup} \ d(A,b)   \right\rbrace.$$

Then $H$ is an extended quasi-pseudometric on $\mathscr{P}_0(X)$. Moreover, we know from \cite{ah} that, the restriction of $H$ to $S_{cl}(X)= \{A\subseteq X: A= cl_{\tau (d)}A \cap cl_{\tau (d^{-1})}A  \}$ is an extended $T_0$-quasi- pseudometric. We shall denote by $CB(X)$ the collection of all nonempty bounded and $\tau(d)$-closed subsets of $X$.

We complete this section by the following lemma.
\begin{lem} (\cite{rico})
Let $(X,d)$ be a quasi-pseudometric space. For every fixed $x\in X,$ the mapping $y\mapsto d(x,y)$ is $\tau(d)$-usc and $\tau(d^{-1})$-lsc. For every fixed $y\in X,$ the mapping $x\mapsto d(x,y)$ is $\tau(d)$-lsc and  $\tau(d^{-1})$-usc.

\end{lem}

For the convenience of the reader, we recall the following.
\begin{dfn}(\cite{rico})
Let $F:X\to 2^X$ be a set-valued map. An element $x\in X$ is said to be 
\begin{enumerate}
\item[(i)] a fixed point of $F$ if $x\in Fx$,
\item[(ii)] a startpoint of $F$ if $H(\{x\},Fx)=0$,
\item[(iii)] an endpoint of $F$ if $H(Fx,\{x\})=0$.
\end{enumerate}
\end{dfn}

We also give here the first result of the theory as they appear in the original paper \cite{rico}.

\begin{theor}\label{result2}(Gaba \cite{rico})
Let $(X,d)$ be a left $K$-complete quasi-pseudometric space. Let $F:X\to CB(X)$ be a set-valued map and $f:X\to \mathbb{R}$ as $f(x)=H(\{x\},Fx)$. If there exists $c\in (0,1)$ such that for all $x\in X$ there exists $y\in Fx$ satisfying

\begin{equation}
H(\{y\},Fy) \leq c (d(x,y)),
\end{equation}
then $T$ has a startpoint.
\end{theor}

Motivated by this result, Gaba proved the following theorems in \cite{ricoo}

\begin{theor}\label{result3}(Gaba \cite{ricoo})
Let $(X,d)$ be a left $K$-complete quasi-pseudometric space. Let $T:X\to CB(X)$ be a set-valued map and define $f:X\to \mathbb{R}$ as $f(x)=H(\{x\},Tx)$. Let $\Phi: [0,\infty)\to [0,1)$ a function such that $\underset{r\to t^+}{\limsup} \ \Phi(r)<1$ for each $t\in [0,\infty)$. Moreover, assume that for any $x\in X$ there exists $y\in Tx$ satisfying

\begin{equation}
H(\{y\},Ty) \leq \Phi(H(\{x\},\{y\}))H(\{x\},\{y\}),
\end{equation}
then $T$ has a startpoint.
\end{theor}

\begin{theor}\label{result4}(Gaba \cite{ricoo})
Let $(X,d)$ be a left $K$-complete quasi-pseudometric space. Let $T:X\to CB(X)$ be a set-valued map and define $f:X\to \mathbb{R}$ as $f(x)=H(\{x\},Tx)$. Let $b:[0,\infty)\to [a,1),$ $a>0$ be a non-decreasing function. Let $\Phi: [0,\infty)\to [0,1)$ be a function such that  $\Phi(t)<b(t)$ for each $t\in [0,\infty)$ and $\underset{r\to t^+}{\limsup}\ \Phi(r)< \underset{r\to t^+}{\limsup}\ b(r)$ for each $t\in [0,\infty)$. Moreover, assume that for any $x\in X$ there exists $y\in Tx$ satisfying

\begin{equation}
H(\{y\},Ty) \leq \Phi(H(\{x\},\{y\}))H(\{x\},\{y\}),
\end{equation}
then $T$ has a startpoint.
\end{theor}

We shall denote by $B(X)$ the collection of all nonempty bounded subsets of $X$.

\section{Main results}

We can state our first result. 

\begin{theor}\label{Result1}
Let $(X,d)$ be a left $K$-complete quasi-pseudometric space. Let $T:X\to B(X)$ be a set-valued map and define $f:X\to \mathbb{R}$ as $f(x)=H(\{x\},Tx)$. Suppose that there exist functions $\varphi:[0,\infty)\to [0,1), \eta:[0,\infty)\to [b,1),$ $0<b<1$ such that 

\[
\varphi(t) < \eta(t),\ \  \underset{r\to t^+}{\limsup}\ \frac{\varphi(r)}{\eta(r)} <1 \ \ \ \text{ for each } t\in [0,\infty).
\]
Moreover, assume that for any $x\in X$ there exists $y\in Tx$ satisfying

\begin{equation}\label{eq1}
f(y) \leq \varphi(f(x)) H(\{x\},\{y\}),
\end{equation}
then $T$ has a startpoint.
\end{theor}

\begin{pf}

Observe that for any $x\in X$ and $y\in Tx$, we have 
\begin{equation}\label{eq2}
\eta(f(x)) H(\{x\},\{y\}) \leq f(x).
\end{equation}

The proof follows similar patterns as in Theorem \ref{result4}. Let $x_0 \in X$, there exists $x_1 \in Tx_0\subseteq X$ such that

\[
\eta(f(x_0)) H(\{x_0\},\{x_1\}) \leq f(x_0) \text{ and } f(x_1) \leq \varphi(f(x_0)) H(\{x_0\},\{x_1\}).
\]

Then

\[
f(x_1) \leq \frac{\varphi(f(x_0))}{\eta(f(x_0))}\eta(f(x_0)) H(\{x_0\},\{x_1\}) \leq  \frac{\varphi(f(x_0))}{\eta(f(x_0))}f(x_0).
\]

In this manner we can build the sequence $(x_n) \subseteq X$ such that for $n\in \mathbb{N},$ $x_{n+1} \in Tx_n\subseteq X$, with

\begin{equation}\label{eq3}
\eta(f(x_n)) H(\{x_n\},\{x_{n+1}\}) \leq f(x_n).
\end{equation}
and

\begin{equation}\label{eq4}
f(x_{n+1}) \leq \frac{\varphi(f(x_n))}{\eta(f(x_n))}f(x_n).
\end{equation}
From \eqref{eq4}, it follows that $(f(x_n))$ is a decreasing sequence of positive real numbers, so there exists $\delta\geq 0$ such that $f(x_n) \to \delta$. If we let $\beta = \underset{n\to \infty}{\limsup}\ \frac{\varphi(f(x_n))}{\eta(f(x_n))}<1$, then for $q=\frac{\beta+1}{2}<1$, there is $n_0 \in \mathbb{N}$ such that 

\[
\frac{\varphi(f(x_n))}{\eta(f(x_n))}<q \  \ \ \  \text{for all } n \geq n_0.
\]

Thus
\begin{equation}\label{eq5}
f(x_{n+1}) \leq q^{n-n_0+1}f(x_n), \  \ \ \ \text{for all } n \geq n_0,
\end{equation}
and so 
\[
H(\{x_n\},\{x_{n+1}\}) \leq \frac{1}{\eta(f(x_n))}f(x_n)\leq \frac{1}{b	}q^{n-n_0}f(x_{n_0}) \ \text{for all } n \geq n_0.
\]
Hence $(x_n)$ is a left $K$-Cauchy sequence. According to the left $K$-completeness of $(X,d),$ there exists $x^* \in X$ such that $x_n \overset{d}{\longrightarrow} x^*$.

Observe that letting $n\to \infty $ in \eqref{eq5}, we get that the sequence $(fx_n)=(H(\{x_{n}\},Tx_{n}))$ converges to $0$. Since $f$ is $\tau(d)$-lower semicontinuous (as supremum of $\tau(d)$-lower semicontinuous functions), we have 

 \[
 0\leq f(x^*) \leq \underset{n\to \infty}{\liminf} f(x_n)=0.
 \]
Hence $f(x^*)=0$, i.e. $H(\{x^*\},Tx^*)=0$.

\noindent
This completes the proof.

\end{pf}

Next we give an example of mapping which satisfies the hypotheses of Theorem \ref{Result1} but does not fulfil the conditions of Theorem \ref{result3}. This example actually motivates our extension.

\begin{ex}
Let $X=[0,10]$ and $d:X\times X \to \mathbb{R}$ be the mapping defined by $d(a,b)=\max\{a-b,0\}$. Then $d$ is a $T_0$-quasi-pseudometric on $X$. Moreover, since any sequence in  $(X,d)$ $d$-converges to $0$, $(X,d)$ is left $K$-complete. Let $T:X\to B(X)$ be such that

\[T(x) = 
  \begin{cases}
   \left\lbrace 4,5\right\rbrace, & \text{if } x=6 ,  \\
    \left\lbrace \frac{x	}{2} \right\rbrace, &\text{if } x \in [0,10]\setminus \{6\}.
 \end{cases}
\]

An explicit computation of $f(x)=H(\{x\},Tx)$ gives

\[f(x) = 
  \begin{cases}
   2 , & \text{if } x=6 ,   \\
     \frac{x}{2}, &\text{if } x \in [0,10]\setminus \{6\}.
 \end{cases}
\]
In view of Theorem \eqref{result3}, we set $\Phi(x)= \frac{1}{2}$ for every $x\geq 0$. So, for $x=6$, $T6=\{4,5\}$  and $f(6)=2$. 

If $y=4$ we have $\Phi(H(\{x\},\{y\}))=\Phi(d(6,4)) = \Phi(2)=\frac{1}{2}$ and

\[H(\{y\},Ty) = d(4,2)=2 > \Phi(H(\{x\},\{y\}))H(\{x\},\{y\})= 1.\]

If $y=5$ we have $\Phi(H(\{x\},\{y\}))=\Phi(d(6,5)) = \Phi(1)=\frac{1}{2}$ and 
\[H(\{y\},Ty) = d\left(5,\frac{5}{2}\right)=\frac{5}{2}> \Phi(H(\{x\},\{y\}))H(\{x\},\{y\})= \frac{1}{2}.\]

Hence the Theorem \eqref{result3} cannot be applied for $T$.

However, it is a simple exercise to show that for $\varphi(x)= \frac{1}{2}$ and $\eta(x)= \frac{2}{3}$ for every $x\geq 0$, the mapping $T$ satisfies the hypotheses of Theorem \eqref{Result1}. Hence $T$ has a starpoint which is $0.$

\end{ex}
We have the logical corollaries.

\begin{cor}
Let $(X,d)$ be a right $K$-complete quasi-pseudometric space. Let $T:X\to B(X)$ be a set-valued map and define $f:X\to \mathbb{R}$ as $f(x)=H(Tx,\{x\})$. Suppose that there exist functions $\varphi:[0,\infty)\to [0,1), \eta:[0,\infty)\to [b,1),$ $0<b<1$ such that 

\[
\varphi(t) < \eta(t),\ \underset{r\to t^+}{\limsup}\ \frac{\varphi(r)}{\eta(r)} <1 \ \ \text{ for each } t\in [0,\infty).
\]
Moreover, assume that for any $x\in X$ there exists $y\in Tx$ satisfying

\begin{equation}
f(y) \leq \varphi(f(x)) H(\{y\},\{x\}),
\end{equation}
then $T$ has an endpoint.
\end{cor}

\begin{cor}
Let $(X,d)$ be a bicomplete quasi-pseudometric space. Let $T:X\to B(X)$ be a set-valued map and define $f:X\to \mathbb{R}$ as $f(x)=H^s(Tx,\{x\})$. Suppose that there exist functions $\varphi:[0,\infty)\to [0,1), \eta:[0,\infty)\to [b,1),$ $0<b<1$ such that 

\[
\varphi(t) < \eta(t),\ \ \underset{r\to t^+}{\limsup}\ \frac{\varphi(r)}{\eta(r)} <1  \ \ \ \text{ for each } t\in [0,\infty).
\]
Moreover, assume that for any $x\in X$ there exists $y\in Tx$ satisfying

\begin{equation}\label{sure}
f(y) \leq  \min \{\varphi(H(\{x\},Tx))H(\{x\},\{y\}), \varphi(H(Tx,\{x\}))H(\{y\},\{x\}) \},
\end{equation}
then $T$ has a fixed point.
\end{cor}

\begin{pf}
We give here the main idea of the proof. Observe that the inequality \eqref{sure} guarantees that the sequence $(x_n)$ constructed in the proof of Theorem \ref{Result1} is a $d^s$-Cauchy sequence and hence $d^s$-converges to some $x^*$. Using the fact that $f$ is $\tau(d^s)$-lower semicontinuous (as supremum of $\tau(d^s)$-continuous functions), we have 

 \[
 0\leq f(x^*) \leq \underset{n\to \infty}{\liminf} f(x_n)=0.
 \]
Hence $f(x^*)=0$, i.e. $H(\{x^*\},Tx^*)=0=H(Tx^*,\{x^*\})$, and we are done.
\end{pf}

The Theorem \ref{result4} imposes a monotonicity condition on the bounding function $b$. We rewrite the original assumptions to obtain the following result.

\begin{theor}\label{Result2}
Let $(X,d)$ be a left $K$-complete quasi-pseudometric space. Let $T:X\to B(X)$ be a set-valued map and define $f:X\to \mathbb{R}$ as $f(x)=H(\{x\},Tx)$. Suppose that there exist functions $\varphi:[0,\infty)\to [0,1), \eta:[0,\infty)\to [b,1],$ $0<b<1$ such that 

\begin{equation}\label{eqs9}
\varphi(t) < \eta(t) \ \ \ \ \text{ for each } t\in [0,\infty),
\end{equation}

\begin{equation}\label{eqs10}
\underset{r\to t^+}{\limsup}\ \frac{\varphi(r)}{\eta(r)} <1 \ \ \ \ \text{ for each } t\in [0,\infty).
\end{equation}

Moreover, assume that for any $x\in X$ there exists $y\in Tx$ satisfying

\begin{equation}\label{eqs12}
f(y) \leq \varphi(H(\{x\},\{y\})) H(\{x\},\{y\}),
\end{equation}
then $T$ has a startpoint.
\end{theor}

\begin{pf}
Observe that for any $x\in X$ and $y\in Tx$, we have 
\begin{equation}\label{eq11}
\eta(H(\{x\},\{y\})) H(\{x\},\{y\}) \leq f(x).
\end{equation}

As in the proof of Theorem \ref{Result1}, we can build a sequence $(x_n) \subseteq X$ uch that for $n\in \mathbb{N},$ $x_{n+1} \in Tx_n\subseteq X$, with
\begin{equation}\label{eqs13}
\eta(H(\{x_n\},\{x_{n+1}\})) H(\{x_n\},\{x_{n+1}\}) \leq f(x_n),
\end{equation}
and 
\begin{equation}\label{eqs14}
f(x_{n+1}) \leq \frac{\varphi(H(\{x_n\},\{x_{n+1}\}))}{\eta(H(\{x_n\},\{x_{n+1}\}))}f(x_n).
\end{equation}

%%%%%%

Set $\Phi(t)=\frac{\varphi}{\eta}(t)$ and since $\Phi(H(\{x\},\{y\})) <1$ for any $x,y \in X$, it follows that $2-\Phi(H(\{x\},\{y\}))>1$ for any $x,y \in X$. Hence 

\begin{equation}\label{conds2}
H(\{x\},\{y\}) \leq [2-\Phi(H(\{x\},\{y\}))] H(\{x\},Tx),
\end{equation}

for any $x\in X$ and $y \in Tx$.

This entails that for the sequence $(x_n)$, we have

\begin{equation}\label{cond1}
H(\{x_n\},\{x_{n+1}\}) \leq [2-\Phi(H(\{x_n\},\{x_{n+1}\}))]H(\{x_n\},Tx_n),
\end{equation}
and 
\begin{equation}\label{cond2}
H(\{x_{n+1}\},Tx_{n+1}) \leq \Psi(H(\{x_n\},\{x_{n+1}\}))H(\{x_n\},Tx_n), \ n=1,2,\cdots
\end{equation}
where $\Psi:[0,\infty) \to [0,1)$ is defined by $$\Psi(t)=\Phi(t)(2-\Phi(t)).$$

For simplicity, denote $d_n:=H(\{x_n\},\{x_{n+1}\})$ and $D_n:=H(\{x_n\},Tx_n)$ for all $n\geq 0.$
So from \eqref{cond2} we can write
\begin{equation}\label{majo}
D_{n+1}\leq \Psi(d_n)D_n\leq D_n
\end{equation}
for all $n\geq 0$. Hence $(D_n)$ is a strictly decreasing sequence, hence there exists $\delta\geq 0$ such that 
\begin{equation}\label{limdn}
\underset{n\to \infty}{\lim} D_n =\delta.
\end{equation}

From \eqref{cond1}, it is easy to see that 
\begin{equation}\label{limdn1}
d_n < 2D_n.
\end{equation}
Thus the sequence $(d_n)$ is bounded and so there is $d\geq 0$ such that $\underset{n\to \infty}{\limsup} \ d_n =d,$ and hence a subsequence $(d_{n_{k}})$ of $(d_n)$ such that $\underset{k\to \infty}{\lim} d_{n_{k}} =d^+$. From \eqref{cond2} we have $D_{n_k+1}\leq \Psi(d_{n_k})D_{n_k}$ and thus

\begin{align*}
\delta = \underset{k\to \infty}{\limsup}\ D_{n_k+1} & \leq (\underset{k\to \infty}{\limsup} \ \Psi(d_{n_k})) (\underset{k\to \infty}{\limsup} \ D_{n_k} ) \\
             & \leq \underset{d_{n_k}\to d^+}{\limsup}\ \Psi(d_{n_k}) \delta .
\end{align*}
This together with the fact $\underset{r\to t^+}{\limsup}\ \Phi(r)<1$ for each $t\in [0,\infty)$ imply that $\delta=0$. Then from \eqref{limdn} and \eqref{limdn1} we derive that $\underset{n\to \infty}{\lim} d_n =0$.

\vspace*{0.3cm}

\underline{Claim 1} $(x_n)$ is a left $K$-Cauchy sequence.

Now let $\alpha := \underset{d_{n}\to 0^+}{\limsup}\ \Psi(d_{n})$ and $q$ such that $\alpha<q<1$. This choice of $q$ is always possible since $\alpha<1$. Then there is $n_0$ such that $\Psi(d_n)<q $ for all $n\geq n_0$. So from \eqref{majo} we have $D_{n+1}\leq q D_n$ for all $n\geq n_0$. Then by induction we get $D_n \leq q^{n-n_0}D_{n_0}$ for all $n\geq n_0+1$. Combining this and the inequality \eqref{limdn1} we get

\[
\sum_{k=n_0}^{m}d(x_k,x_{k+1})\leq 2 \sum_{k=n_0}^{m} q^{k-n_0}D_{n_0} \leq 2 \frac{1}{1-q}D_{n_0}
\]
for all $m>n\geq n_0+1$. Hence $(x_n)$ is a left $K$-Cauchy sequence.

According to the left $K$-completeness of $(X,d),$ there exists $x^* \in X$ such that $x_n \overset{d}{\longrightarrow} x^*$.

\vspace*{0.3cm}

\underline{Claim 2} $x^*$ is a startpoint of $T$.

Observe that the sequence $D_n=(fx_n)=(H(\{x_{n}\},Tx_{n}))$ converges to $0$. Since $f$ is $\tau(d)$-lower semicontinuous (as supremum of $\tau(d)$-lower semicontinuous functions), we have 

 \[
 0\leq f(x^*) \leq \underset{n\to \infty}{\liminf} f(x_n)=0.
 \]
Hence $f(x^*)=0$, i.e. $H(\{x^*\},Tx^*)=0$.

This completes the proof.
\end{pf}

\begin{cor}
Let $(X,d)$ be a right $K$-complete quasi-pseudometric space. Let $T:X\to B(X)$ be a set-valued map and define $f:X\to \mathbb{R}$ as $f(x)=H(Tx,\{x\})$. Suppose that there exist functions $\varphi:[0,\infty)\to [0,1), \eta:[0,\infty)\to [b,1],$ $0<b<1$ such that 

\begin{equation}
\varphi(t) < \eta(t) \ \ \ \ \text{ for each } t\in [0,\infty),
\end{equation}

\begin{equation}
\underset{r\to t^+}{\limsup}\ \frac{\varphi(r)}{\eta(r)} <1 \ \ \ \ \text{ for each } t\in [0,\infty).
\end{equation}

Moreover, assume that for any $x\in X$ there exists $y\in Tx$ satisfying

\begin{equation}
f(y) \leq \varphi(H(\{y\},\{x\})) H(,\{y\},\{x\}),
\end{equation}
then $T$ has an endpoint.
\end{cor}

\begin{cor}
Let $(X,d)$ be a right $K$-complete quasi-pseudometric space. Let $T:X\to B(X)$ be a set-valued map and define $f:X\to \mathbb{R}$ as $f(x)=H^s(Tx,\{x\})$. Suppose that there exist functions $\varphi:[0,\infty)\to [0,1), \eta:[0,\infty)\to [b,1],$ $0<b<1$ such that 

\begin{equation}
\varphi(t) < \eta(t) \ \ \ \ \text{ for each } t\in [0,\infty),
\end{equation}

\begin{equation}
\underset{r\to t^+}{\limsup}\ \frac{\varphi(r)}{\eta(r)} <1 \ \ \ \ \text{ for each } t\in [0,\infty).
\end{equation}

Moreover, assume that for any $x\in X$ there exists $y\in Tx$ satisfying

\begin{equation}
f(y) \leq \min \{\varphi(a)a,\varphi(b)b\} ,
\end{equation}
where $a=H(\{y\},\{x\})$ and $b= H(\{y\},\{x\}). $
Then $T$ has a fixed point.
\end{cor}

Another approach worth investigating would be to try to generalize even further conditions \eqref{eq1} and \eqref{eqs12}. In this direction, we state the new results.

\begin{theor}\label{Result3}
Let $(X,d)$ be a left $K$-complete quasi-pseudometric space. Let $T:X\to B(X)$ be a set-valued map and define $f:X\to \mathbb{R}$ as $f(x)=H(\{x\},Tx)$. Suppose that there exist functions $\varphi:[0,\infty)\to [0,\infty), \eta:[0,\infty)\to [0,\infty)$ such that 

\[ \varphi \text{ is non-decreasing
 and} \ 
\varphi(t) < \eta(t), \ \underset{r\to t^+}{\limsup}\ \frac{\varphi(r)}{\eta(r)} <1 \ \text{ for each } t\in [0,\infty).
\]
Moreover, assume that for any $x\in X$ there exists $y\in Tx$ satisfying

\begin{equation}\label{eqs18}
\eta(H(\{x\},\{y\})) \leq f(x),
\end{equation}
and 
\begin{equation}\label{eqs19}
f(y) \leq \varphi(f(x)) .
\end{equation}

Then $T$ has a startpoint.
\end{theor}

\begin{pf}
Let $x_0 \in X$. We can choose $x_1\in Tx_0$ such that

\[
\eta(H(\{x_0\},\{x_1\})) \leq f(x_0)  \text { and } f(x_1) \leq \varphi(f(x_0)) .
\]
In this way, we build the sequence $(x_n) \subseteq X$ such that for $n\in \mathbb{N},$ $x_{n+1} \in Tx_n$, with

\begin{equation}\label{eqs20}
\eta(H(\{x_n\},\{x_{n+1}\})) \leq f(x_n) 
\end{equation}
and 
\begin{equation}\label{eqs21}
f(x_{n+1}) \leq \varphi(f(x_n)).
\end{equation}

\end{pf}

Since 
\begin{align*}
\varphi(H(\{x_{n+1}\},\{x_{n+2}\})) & < \eta(H(\{x_{n+1}\},\{x_{n+2}\})) \\
               & \leq f(x_{n+1}) \leq \varphi(f(x_{n})) \\
                & \leq \varphi(H(\{x_{n}\},\{x_{n+1}\})),         
\end{align*}

it follows that $\varphi(H(\{x_{n+1}\},\{x_{n+2}\})) < \varphi(H(\{x_{n}\},\{x_{n+1}\}))$. Because $\varphi$ is non-decreasing, the sequence $(H(\{x_{n}\},\{x_{n+1}\}))$ is decreasing. Because  the sequence $(H(\{x_{n}\},\{x_{n+1}\}))$ is bounded from below, it converges. Using \eqref{eqs20} and \eqref{eqs21} we have that 

\[
f(x_{n+1}) \leq \frac{\varphi(H(\{x_n\},\{x_{n+1}\}))}{\eta(H(\{x_n\},\{x_{n+1}\}))}f(x_n).
\]
Due to $$ \underset{r\to t^+}{\limsup}\ \frac{\varphi(r)}{\eta(r)} <1, $$ there exist $q\in (0,1)$ and $n_0\in \mathbb{N}$

\[
\frac{\varphi(H(\{x_n\},\{x_{n+1}\}))}{\eta(H(\{x_n\},\{x_{n+1}\}))} <q \ \ \ \text{ for all } n\geq n_0 .
\]

Thus 
\[
f(x_{n+1}) \leq q^{n-n_0+1}f(x_{n_0}) \ \ \ \text{ for all } n\geq n_0.
\]

For $ n\geq n_0+1,$
\begin{align*}
\varphi(H(\{x_{n}\},\{x_{n+1}\})) & < \eta(H(\{x_{n}\},\{x_{n+1}\})) \\
               & \leq f(x_{n}) \leq \varphi(f(x_{n-1})) \\
                & \leq \varphi(q^{n-n_0-1}f(x_{n_0})).         
\end{align*}

Because $\varphi$ is non-decreasing, $ H(\{x_{n}\},\{x_{n+1}\}) \leq q^{n-n_0-1}f(x_{n_0}) $. It is therefore easy to see that $(x_n)$  is a left $K$-Cauchy sequence and its limit is a startpoint for $T$.

\begin{cor}
Let $(X,d)$ be a right $K$-complete quasi-pseudometric space. Let $T:X\to B(X)$ be a set-valued map and define $f:X\to \mathbb{R}$ as $f(x)=H(Tx,\{x\})$. Suppose that there exist functions $\varphi:[0,\infty)\to [0,\infty), \eta:[0,\infty)\to [0,\infty)$ such that 

\[ \varphi \text{ is non-decreasing 
and} \ 
\varphi(t) < \eta(t), \ \underset{r\to t^+}{\limsup}\ \frac{\varphi(r)}{\eta(r)} <1  \ \text{ for each } t\in [0,\infty).
\]
Moreover, assume that for any $x\in X$ there exists $y\in Tx$ satisfying

\begin{equation}
\eta(H(\{y\},\{x\})) \leq f(x),
\end{equation}
and 
\begin{equation}
f(y) \leq \varphi(f(x)) .
\end{equation}

Then $T$ has an endpoint.
\end{cor}

\begin{cor}
Let $(X,d)$ be a bicomplete quasi-pseudometric space. Let $T:X\to B(X)$ be a set-valued map and define $f:X\to \mathbb{R}$ as $f(x)=H^s(\{x\},Tx)$. Suppose that there exist functions $\varphi:[0,\infty)\to [0,\infty), \eta:[0,\infty)\to [0,\infty)$ such that 

\[ \varphi \text{ is non-decreasing and} \ 
\varphi(t) < \eta(t),\ \underset{r\to t^+}{\limsup}\ \frac{\varphi(r)}{\eta(r)} <1  \ \text{ for each } t\in [0,\infty).
\]
Moreover, assume that for any $x\in X$ there exists $y\in Tx$ satisfying

\begin{equation}
\eta(H^s(\{x\},\{y\}))
\leq \min \{ H(\{x\},Tx), H(Tx,\{x\})\},
\end{equation}
and 
\begin{equation}
f(y) \leq \min \{\varphi(H(\{x\},Tx)), \varphi(H(Tx,\{x\}))\}.
\end{equation}

Then $T$ has a fixed point.
\end{cor}

\begin{theor}\label{Result4}
Let $(X,d)$ be a left $K$-complete quasi-pseudometric space. Let $T:X\to B(X)$ be a set-valued map and define $f:X\to \mathbb{R}$ as $f(x)=H(\{x\},Tx)$. Suppose that there exist functions $\varphi:[0,\infty)\to [0,\infty), \eta:[0,\infty)\to [0,\infty)$ such that 

\[ \eta \text{ is non-decreasing and} \ 
\varphi(t) < \eta(t), \ \underset{r\to t^+}{\limsup}\ \frac{\varphi(r)}{\eta(r)} <1 \ \text{ for each } t\in [0,\infty).
\]
Moreover, assume that for any $x\in X$ there exists $y\in Tx$ satisfying

\begin{equation}\label{eqs18}
\eta(H(\{x\},\{y\})) \leq f(x),
\end{equation}
and 
\begin{equation}\label{eqs19}
f(y) \leq \varphi(f(x)) .
\end{equation}

Then $T$ has a startpoint.
\end{theor}
\begin{pf}
We build the sequence $(x_n)\subseteq X$ as in the proof of Theorem \ref{Result3}.
Since $\eta$ is non-decreasing we obtain we obtain that for $n\in \mathbb{N}$,

\[
f(x_{n+1}) \leq \frac{\varphi(f(x_{n}))}{\eta(f(x_{n}))}f(x_n).
\] 
Hence $(f(x_n))$ is decreasing. Because  the sequence $(f(x_n))$ is bounded from below, it converges. Again there exist $q\in (0,1)$ and $n_0\in \mathbb{N}$ such that

\[
f(x_{n+1}) \leq q^{n-n_0+1}f(x_{n_0}) \ \ \ \text{ for all } n\geq n_0.
\]

For $ n\geq n_0+1,$
\begin{align*}
\eta(H(\{x_{n}\},\{x_{n+1}\})) & \leq f(x_{n})  \leq \varphi(f(x_{n-1}))  \\
               & \leq \eta(f(x_{n-1}))  \\
                & \leq \eta(q^{n-n_0-1}f(x_{n_0})).         
\end{align*}
But $\eta$ is non-decreasing, so $H(\{x_{n}\},\{x_{n+1}\}) \leq q^{n-n_0-1}f(x_{n_0})$. As above we can show that $(x_n)$  is a left $K$-Cauchy sequence and its limit is a startpoint for $T$.

\end{pf}

\begin{cor}
Let $(X,d)$ be a right $K$-complete quasi-pseudometric space. Let $T:X\to B(X)$ be a set-valued map and define $f:X\to \mathbb{R}$ as $f(x)=H(Tx,\{x\})$. Suppose that there exist functions $\varphi:[0,\infty)\to [0,\infty), \eta:[0,\infty)\to [0,\infty)$ such that 

\[ \eta \text{ is non-decreasing and} \ 
\varphi(t) < \eta(t),\ \underset{r\to t^+}{\limsup}\ \frac{\varphi(r)}{\eta(r)} <1 \ \text{ for each } t\in [0,\infty).
\]
Moreover, assume that for any $x\in X$ there exists $y\in Tx$ satisfying

\begin{equation}
\eta(H(\{y\},\{x\})) \leq f(x),
\end{equation}
and 
\begin{equation}
f(y) \leq \varphi(f(x)) .
\end{equation}

Then $T$ has an endpoint.
\end{cor}

\begin{cor}
Let $(X,d)$ be a bicomplete quasi-pseudometric space. Let $T:X\to B(X)$ be a set-valued map and define $f:X\to \mathbb{R}$ as $f(x)=H^s(Tx,\{x\})$. Suppose that there exist functions $\varphi:[0,\infty)\to [0,\infty), \eta:[0,\infty)\to [0,\infty)$ such that 

\[ \eta \text{ is non-decreasing and} \ 
\varphi(t) < \eta(t),\ \underset{r\to t^+}{\limsup}\ \frac{\varphi(r)}{\eta(r)} <1  \ \text{ for each } t\in [0,\infty).
\]
Moreover, assume that for any $x\in X$ there exists $y\in Tx$ satisfying

\begin{equation}
\eta(H^s(\{x\},\{y\})) \leq 
\min \{ H(\{x\},Tx), H(Tx,\{x\})\},
\end{equation}
and 
\begin{equation}
f(y) \leq \min \{ \varphi(H(Tx,\{x\})), \varphi(H(\{x\},Tx))\} .
\end{equation}
Then $T$ has a fixed point.
\end{cor}

In the sequel, we prove two related theorems.

\begin{theor}\label{Result5}
Let $(X,d)$ be a left $K$-complete quasi-pseudometric space. Let $T:X\to B(X)$ be a set-valued map and define $f:X\to \mathbb{R}$ as $f(x)=H(\{x\},Tx)$. Suppose that there exist functions $\varphi:[0,\infty)\to [0,\infty), \eta:[0,\infty)\to [0,\infty)$ such that $\varphi \text{ is continuous, non-decreasing},$ 

\[  \
\varphi(t) < \eta(t) \leq t, \ \ \underset{r\to t^+}{\limsup}\ \frac{\varphi(r)}{\eta(r)} <1 \ \ \ \text{ for each } t\in [0,\infty).
\]
Moreover, assume that for any $x\in X$ there exists $y\in Tx$ satisfying

\begin{equation}
\eta(H(\{x\},\{y\})) \leq f(x),
\end{equation}
and 
\begin{equation}
f(y) \leq \varphi(f(x)) .
\end{equation}
Then $T$ has a startpoint.
\end{theor}

\begin{pf}
Again we build the sequence $(x_n)$ with $x_{n+1} \in Tx_n$ such that \eqref{eqs20} and \eqref{eqs21} hold. Without loss of generality, we can always assume that $H(\{x_{n}\},\{x_{n+1}\})>0$ and 
$f(x_n)>0$ for $n\in \mathbb{N},$ (because otherwise, the proof is already complete).

Now, let $t>0$, since $0\leq \varphi(t)<t$, the sequence $(\varphi^n(t))$ is decreasing and bounded from below by $0$. Suppose its limit is $\zeta >0$, then 
\[
\zeta = \underset{n\to \infty}
{\lim}\ \varphi^n(t)= \varphi 
\left(\underset{n\to \infty}{\lim}\ 
\varphi^{n-1}(t)\right) = \varphi(\zeta) < \zeta,
\]
which is a contradiction. Therefore $\underset{n\to \infty}{\lim}\ \varphi^n(t)=0 \text{ for each } t\in [0,\infty).$

From \eqref{eqs21}, we have that 
\[
f(x_{n+1}) \leq \varphi(f(x_n)) \leq \cdots \leq \varphi^{n+1}(f(x_0)),
\]
and $\underset{n\to \infty}{\lim}\ f(x_n)=0$.

Since $\varphi$ is non-decreasing
\begin{align*}
\varphi(H(\{x_{n+1}\},\{x_{n+2}\})) & < \eta(H(\{x_{n+1}\},\{x_{n+2}\})) \\
               & \leq f(x_{n+1}) \leq \varphi(f(x_{n})) \\
                & \leq \varphi(H(\{x_{n}\},\{x_{n+1}\})),         
\end{align*}

the sequence $(H(\{x_{n+1}\},\{x_{n+2}\}))$ is decreasing. If $\underset{n\to \infty}{\lim}\ H(\{x_{n}\},\{x_{n+1}\})=\xi >0$, we have, by letting $n\to \infty$ in $\varphi(H(\{x_{n}\},\{x_{n+1}\})) < f(x_n)$ that $\varphi(\xi)=0$. On the other hand, since $\underset{n\to \infty}{\lim}\ f(x_n)=0,$ there exists $n_1 \in \mathbb{N}$ such that $f(x_{n_1}) < \xi$ and hence $\varphi(f(x_{n_1}))=0$. This means that $f(x_{n_1+1})=0$. In this way we obtain a startpoint. Therefore, we may consider that $\underset{n\to \infty}{\lim}\ H(\{x_{n}\},\{x_{n+1}\})=0$. Continuing as in the proof of 
Theorem \eqref{Result3} one can show that $T$ is not startpoint free.
\end{pf}

\begin{cor}
Let $(X,d)$ be a right $K$-complete quasi-pseudometric space. Let $T:X\to B(X)$ be a set-valued map and define $f:X\to \mathbb{R}$ as $f(x)=H(Tx,\{x\})$. Suppose that there exist functions $\varphi:[0,\infty)\to [0,\infty), \eta:[0,\infty)\to [0,\infty)$ such that $\varphi \text{ is continuous, non-decreasing},$ 

\[  \
\varphi(t) < \eta(t) \leq t, \ \ \underset{r\to t^+}{\limsup}\ \frac{\varphi(r)}{\eta(r)} <1 \ \ \ \text{ for each } t\in [0,\infty).
\]
Moreover, assume that for any $x\in X$ there exists $y\in Tx$ satisfying

\begin{equation}
\eta(H(\{y\},\{x\})) \leq f(x),
\end{equation}
and 
\begin{equation}
f(y) \leq \varphi(f(x)) .
\end{equation}
Then $T$ has an endpoint.
\end{cor}

\begin{cor}
Let $(X,d)$ be a bicomplete quasi-pseudometric space. Let $T:X\to B(X)$ be a set-valued map and define $f:X\to \mathbb{R}$ as $f(x)=H^s(\{x\},Tx)$. Suppose that there exist functions $\varphi:[0,\infty)\to [0,\infty), \eta:[0,\infty)\to [0,\infty)$ such that $\varphi \text{ is continuous, non-decreasing},$ 

\[  \
\varphi(t) < \eta(t) \leq t, \ \ \underset{r\to t^+}{\limsup}\ \frac{\varphi(r)}{\eta(r)} <1 \ \ \ \text{ for each } t\in [0,\infty).
\]
Moreover, assume that for any $x\in X$ there exists $y\in Tx$ satisfying

\begin{equation}
\min \{\eta(H(\{y\},\{x\})) 
, \eta(H(\{x\},\{y\}))\}
 \leq \min\{H(\{x\},Tx),H(Tx,\{x\})\},
\end{equation}
and 
\begin{equation}
f(y) \leq \min\{\varphi(H(\{x\},Tx)),\varphi(H(Tx,\{x\}))\}.
\end{equation}
Then $T$ has a fixed point.
\end{cor}

Using a similar argument as above we can prove the following result.

\begin{theor}\label{Result6}
Let $(X,d)$ be a left $K$-complete quasi-pseudometric space. Let $T:X\to B(X)$ be a set-valued map and define $f:X\to \mathbb{R}$ as $f(x)=H(\{x\},Tx)$. Suppose that there exist functions $\varphi:[0,\infty)\to [0,\infty), \eta:[0,\infty)\to [0,\infty)$ such that $\eta \text{ is continuous, non-decreasing},$ 

\[  \
\varphi(t) < \eta(t) < t, \ \ \underset{r\to t^+}{\limsup}\ \frac{\varphi(r)}{\eta(r)} <1 \ \ \ \text{ for each } t\in [0,\infty).
\]
Moreover, assume that for any $x\in X$ there exists $y\in Tx$ satisfying

\begin{equation}
\eta(H(\{x\},\{y\})) \leq f(x),
\end{equation}
and 
\begin{equation}
f(y) \leq \varphi(f(x)) .
\end{equation}
Then $T$ has a startpoint.
\end{theor}

\begin{pf}
We give a sketch of the proof. We use the fact that $\underset{n\to \infty}{\lim}\ \eta^n(t)=0 \text{ for each } t\in [0,\infty)$ and
\begin{align*}
f(x_{n+1}) & \leq \varphi(f(x_n)) < \eta(f(x_n)) \leq \eta ( \varphi(f(x_{n-1}))) \\
           & \leq \eta^2 (f(x_{n-1})) \leq \cdots \leq \eta^{n+1} (f(x_0)), 
\end{align*}
to conclude that $\underset{n\to \infty}{\lim}\ f(x_n)=0.$ As in Theorem \eqref{Result4}, we can prove that $T$ has a startpoint.
\end{pf}

\begin{cor}
Let $(X,d)$ be a right $K$-complete quasi-pseudometric space. Let $T:X\to B(X)$ be a set-valued map and define $f:X\to \mathbb{R}$ as $f(x)=H(Tx,\{x\})$. Suppose that there exist functions $\varphi:[0,\infty)\to [0,\infty), \eta:[0,\infty)\to [0,\infty)$ such that $\eta \text{ is continuous, non-decreasing},$ 

\[  \
\varphi(t) < \eta(t) < t, \ \ \underset{r\to t^+}{\limsup}\ \frac{\varphi(r)}{\eta(r)} <1 \ \ \ \text{ for each } t\in [0,\infty).
\]
Moreover, assume that for any $x\in X$ there exists $y\in Tx$ satisfying

\begin{equation}
\eta(H(\{y\},\{x\})) \leq f(x),
\end{equation}
and 
\begin{equation}
f(y) \leq \varphi(f(x)) .
\end{equation}
Then $T$ has an endpoint.
\end{cor}

\begin{cor}
Let $(X,d)$ be a bicomplete quasi-pseudometric space. Let $T:X\to B(X)$ be a set-valued map and define $f:X\to \mathbb{R}$ as $f(x)=H^s(Tx,\{x\})$. Suppose that there exist functions $\varphi:[0,\infty)\to [0,\infty), \eta:[0,\infty)\to [0,\infty)$ such that $\eta \text{ is continuous, non-decreasing},$ 

\[  \
\varphi(t) < \eta(t) < t, \ \ \underset{r\to t^+}{\limsup}\ \frac{\varphi(r)}{\eta(r)} <1 \ \ \ \text{ for each } t\in [0,\infty).
\]
Moreover, assume that for any $x\in X$ there exists $y\in Tx$ satisfying

\begin{equation}
\eta(H^s(\{x\},\{y\}))
\leq 
\min\{ H(\{x\},Tx),H(Tx,\{x\})\},
\end{equation}
and 
\begin{equation}
f(y) \leq \min\{\varphi(H(Tx,\{x\})),\varphi(H(\{x\},Tx))\} .
\end{equation}
Then $T$ has a fixed point.
\end{cor}

In the above results it would be interesting to replace the condition \eqref{eqs19} by $f(y)\leq \varphi(H(\{x\},\{y\}))$. Pursuing this idea we derive the following two theorems.

\begin{theor}\label{Result7}
Let $(X,d)$ be a left $K$-complete quasi-pseudometric space. Let $T:X\to B(X)$ be a set-valued map and define $f:X\to \mathbb{R}$ as $f(x)=H(\{x\},Tx)$. Suppose that there exist functions $\varphi:[0,\infty)\to [0,\infty), \eta:[0,\infty)\to [0,\infty)$ such that 
$\varphi \text{ is non-decreasing and subadditive}$
\[ \ 
\varphi(t) < \eta(t), \ \underset{r\to t^+}{\limsup}\ \frac{\varphi(r)}{\eta(r)} <1 \ \ \text{ for each } t\in [0,\infty).
\]
Moreover, assume that for any $x\in X$ there exists $y\in Tx$ satisfying

\begin{equation}
\eta(H(\{x\},\{y\})) \leq f(x),
\end{equation}
and 
\begin{equation}
f(y) \leq  \varphi(H(\{x\},\{y\})).
\end{equation}

Then $T$ has a startpoint.
\end{theor}

\begin{pf}
Similarly as before, we can build the sequence $(x_n)$ with $x_{n+1} \in Tx_n$,

\[
\eta(H(\{x_n\},\{x_{n+1}\})) \leq f(x_n) 
\text{ and }
f(x_{n+1}) \leq \varphi(H(\{x_n\},\{x_{n+1}\})).
\]

The sequence $(H(\{x_{n+1}\},\{x_{n+2}\}))$ is decreasing since $\varphi$ is non-decreasing and
\[
\varphi(H(\{x_{n+1}\},\{x_{n+2}\}))  < \eta(H(\{x_{n+1}\},\{x_{n+2}\})) 
               \leq f(x_{n+1})
                 \leq \varphi(H(\{x_{n}\},\{x_{n+1}\})).       
\]
Thus it is convergent. Then there exist $q\in (0,1)$ and $n_0\in \mathbb{N}$ such that

\[
f(x_{n+1}) \leq q^{n-n_0+1}f(x_{n_0}) \ \ \ \text{ for all } n\geq n_0.
\]

For $n\geq n_0$
\[
\varphi(H(\{x_{n}\},\{x_{n+1}\})) < \eta(H(\{x_{n}\},\{x_{n+1}\})) \leq 
f(x_n) \leq q^{n-n_0}f(x_{n_0}) .
\]

For $n\geq n_0$ and $p\in \mathbb{N}$
\begin{align*}
\varphi(H(\{x_{n}\},\{x_{n+p}\})) & \leq \varphi \left( \sum_{k=0}^{p-1} H(\{x_{n+k}\},\{x_{n+k+1}\})  \right) \leq \sum_{k=0}^{p-1} \varphi(H(\{x_{n+k}\},\{x_{n+k+1}\}) ) \\
      & \leq  \sum_{k=0}^{p-1} q^{n-n_0+k}f(x_{n_0}) \leq \frac{q^{n-n_0}}{1-q}f(x_{n_0}).
\end{align*}

Since we may assume that $\varphi(t)>0$ for $t>0$ (otherwise $\varphi(t)=0$) for every $t\geq 0$ and the existence of a startpoint is immediate) we can prove buy contradiction that $(x_n)$ is left $K$-Cauchy sequence and its limit is a startpoint for $T$.
\end{pf}

\begin{cor}
Let $(X,d)$ be a right $K$-complete quasi-pseudometric space. Let $T:X\to B(X)$ be a set-valued map and define $f:X\to \mathbb{R}$ as $f(x)=H(Tx,\{x\})$. Suppose that there exist functions $\varphi:[0,\infty)\to [0,\infty), \eta:[0,\infty)\to [0,\infty)$ such that 
$\varphi \text{ is non-decreasing and subadditive}$
\[ \ 
\varphi(t) < \eta(t), \ \underset{r\to t^+}{\limsup}\ \frac{\varphi(r)}{\eta(r)} <1 \ \text{ for each } t\in [0,\infty).
\]
Moreover, assume that for any $x\in X$ there exists $y\in Tx$ satisfying

\begin{equation}
\eta(H(\{y\},\{x\})) \leq f(x),
\end{equation}
and 
\begin{equation}
f(y) \leq  \varphi(H(\{y\},\{x\})).
\end{equation}

Then $T$ has an endpoint.
\end{cor}

\begin{cor}
Let $(X,d)$ be a bicomplete quasi-pseudometric space. Let $T:X\to B(X)$ be a set-valued map and define $f:X\to \mathbb{R}$ as $f(x)=H^s(\{x\},Tx)$. Suppose that there exist functions $\varphi:[0,\infty)\to [0,\infty), \eta:[0,\infty)\to [0,\infty)$ such that 
$\varphi \text{ is non-decreasing and subadditive}$
\[ \ 
\varphi(t) < \eta(t), \ \underset{r\to t^+}{\limsup}\ \frac{\varphi(r)}{\eta(r)} <1  \ \text{ for each } t\in [0,\infty).
\]
Moreover, assume that for any $x\in X$ there exists $y\in Tx$ satisfying

\begin{equation}
\eta(H^s(\{x\},\{y\}))
\leq \min\{H(\{x\},Tx),H(Tx,\{x\})\} 
,
\end{equation}
and 
\begin{equation}
f(y) \leq \min\{ \varphi(H(\{x\},\{y\})), \varphi(H(\{y\},\{x\})) \}.
\end{equation}

Then $T$ has a fixed point.
\end{cor}

In the same manner we can prove that

\begin{theor}\label{Result8}
Let $(X,d)$ be a left $K$-complete quasi-pseudometric space. Let $T:X\to B(X)$ be a set-valued map and define $f:X\to \mathbb{R}$ as $f(x)=H(\{x\},Tx)$. Suppose that there exist functions $\varphi:[0,\infty)\to [0,\infty), \eta:[0,\infty)\to [0,\infty)$ such that 
$\eta \text{ is non-decreasing and subadditive}$
\[ \ 
\varphi(t) < \eta(t), \ \underset{r\to t^+}{\limsup}\ \frac{\varphi(r)}{\eta(r)} <1 \ \text{ for each } t\in [0,\infty).
\]
Moreover, assume that for any $x\in X$ there exists $y\in Tx$ satisfying

\begin{equation}
\eta(H(\{x\},\{y\})) \leq f(x),
\end{equation}
and 
\begin{equation}
f(y) \leq  \varphi(H(\{x\},\{y\})).
\end{equation}

Then $T$ has a startpoint.
\end{theor}

\begin{cor}
Let $(X,d)$ be a right $K$-complete quasi-pseudometric space. Let $T:X\to B(X)$ be a set-valued map and define $f:X\to \mathbb{R}$ as $f(x)=H(Tx,\{x\})$. Suppose that there exist functions $\varphi:[0,\infty)\to [0,\infty), \eta:[0,\infty)\to [0,\infty)$ such that 
$\eta \text{ is non-decreasing and subadditive}$
\[ \ 
\varphi(t) < \eta(t),\ \underset{r\to t^+}{\limsup}\ \frac{\varphi(r)}{\eta(r)} <1  \ \text{ for each } t\in [0,\infty).
\]
Moreover, assume that for any $x\in X$ there exists $y\in Tx$ satisfying

\begin{equation}
\eta(H(\{y\},\{x\})) \leq f(x),
\end{equation}
and 
\begin{equation}
f(y) \leq  \varphi(H(\{y\},\{x\}).
\end{equation}

Then $T$ has an endpoint.
\end{cor}

\begin{cor}
Let $(X,d)$ be a bicomplete quasi-pseudometric space. Let $T:X\to B(X)$ be a set-valued map and define $f:X\to \mathbb{R}$ as $f(x)=H^s(\{x\},Tx)$. Suppose that there exist functions $\varphi:[0,\infty)\to [0,\infty), \eta:[0,\infty)\to [0,\infty)$ such that 
$\eta \text{ is non-decreasing and subadditive}$
\[ \ 
\varphi(t) < \eta(t), \ \underset{r\to t^+}{\limsup}\ \frac{\varphi(r)}{\eta(r)} <1  \ \text{ for each } t\in [0,\infty).
\]
Moreover, assume that for any $x\in X$ there exists $y\in Tx$ satisfying

\begin{equation}
\eta(H^s(\{x\},\{y\}))
 \leq \min\{  H(\{x\},Tx),
 H(Tx,\{x\})
 \},
\end{equation}
and 
\begin{equation}
f(y) \leq  \min\{\varphi(H(\{x\},\{y\})),\varphi(H(\{y\},\{x\}))\}.
\end{equation}

Then $T$ has a fixed point.
\end{cor}

\begin{rem}
All the results given remain true when we replace accordingly the bicomplete quasi-pseudometric space $(X,d)$ by a left Smyth sequentially complete/left $K$-complete or a right Smyth sequentially complete/right $K$-complete space.
\end{rem}

\end{document}